\documentclass[french, 11pt,a4paper]{article}
\usepackage[utf8]{inputenc}
\usepackage[francais]{babel}
\usepackage[T1]{fontenc}
\usepackage{amsmath}
\usepackage{amsfonts}
\usepackage{amssymb}
\usepackage{amsthm}
\usepackage[left=3cm,right=2cm,top=2cm,bottom=2cm]{geometry}

\newtheorem*{thm}{Théorème}
\newtheorem*{prop}{Proposition}
\newtheorem*{lem}{Lemme}
\newtheorem*{cor}{Corollaire}
\newtheorem*{rmq}{Remarque}

\newcommand{\SL}{\mathrm{SL}}
  \newcommand{\GL}{\mathrm{GL}}

      \newcommand{\SO}{\mathrm{SO}}

\def \dem {\noindent \underline{\sl Démonstration}. }

\begin{document}
\title{Sur les paquets d'Arthur  des groupes classiques et unitaires  non quasi-déployés}
\date{\today}
 \author{Colette Moeglin\\
 CNRS, Institut Mathématique de Jussieu \\
 colette.moeglin@imj-prg.fr
\and
David Renard  \\Centre de Mathématiques
  Laurent Schwartz,  Ecole Polytechnique\\
david.renard@polytechnique.edu}

\maketitle

\begin{abstract} Nous étendons aux groupes orthogonaux et unitaires non quasi-déployés sur un corps local des résultats 
de J. Arthur et de la première auteure établis dans le cas quasi-déployé. 
En particulier, nous obtenons une classification de Langlands complète  pour les représentations tempérées dans le cas $p$-adique.
Nous en déduisons en utilisant l'involution d'Aubert-Schneider-Stuhler un résultat de multiplicité un dans les paquets
unipotents, et par des méthodes globales, le même résultat pour les paquets unipotents dans le cas archimédien.
\medskip

\begin{center} {\bf Abstract}\end{center}  We extend to non quasi-split orthogonal and unitary groups over a local field some results
of J. Arthur and the first author established in the quasi-split case. In  particular, we obtain a full Langlands classification
for tempered representations in the $p$-adic case. Using Aubert-Schneider-Stuhler involution, we deduce from this 
a multiplicity one result for unipotent packets, and by global methods, the same result for unipotent packets in the archimedean case. 

\end{abstract}

\section{Introduction\label{defpaquet}}
Soient $G$ un groupe réductif défini sur un corps local $F$ et $G^*$ sa forme intérieure quasi-déployée.
 En suivant Arthur, on définit les $A$-morphismes de $W'_F\times \SL(2,\mathbb{C})$ dans $^LG$ que l'on note génériquement $\psi$,
  ici $W'_F$ est le groupe de Weil-Deligne de $F$. Dans le cas  $p$-adique qui occupe une bonne place dans ces notes, 
  $W'_F$ est le produit de $W_F$ par $\SL(2,\mathbb{C})$ tandis que si $F$ est un corps archimédien, $W'_F$ est le groupe de Weil de $F$.

Bien que la théorie d'Arthur soit générale, on se limite ici  au cas des groupes orthogonaux et unitaires. Le cas des groupes unitaires est
 d'ailleurs déjà largement fait dans \cite{kalethaandco}. Précisons tout de suite que l'on ne regarde que les formes intérieures pures des groupes
  orthogonaux ou unitaires quasi-déployés. Comme c'est le cas aux places archimédiennes, si on regarde des formes intérieures plus générales,
   le groupe $A(\psi)$ qui intervient ci-dessous n'est pas celui qu'il faut prendre. Ce seul point ne serait pas très grave, mais il y a une autre
    difficulté expliquée par Arthur : c'est l'automorphisme venant du groupe orthogonal  dans le cas des groupes spéciaux orthogonaux pairs 
 qui pose des problèmes. Ainsi le groupe $G$ vient avec le choix d'une forme bilinéaire non dégénérée soit orthogonale soit hermitienne 
 suivant les cas.

 On note $A(\psi)$ le groupe des composantes du centralisateur de $\psi$ dans la composante neutre du $L$-groupe de $G$.
 Ce groupe est commutatif, c'est même un $2$-groupe.
Quand $\psi$ est fixé, on définit $s_\psi$ comme l'image par $\psi$ de l'élément non trivial du centre de
 $\SL(2,\mathbb{C})$. 

On appelle paquet d'Arthur pour $G$ associé à $\psi$ une combinaison linéaire à coefficients complexes de représentations de
 $G\times A(\psi)$,  notée $\pi_G^A(\psi)$ qui vérifie les propriétés énoncées ci-dessous, pour lesquelles il faut quelques notations.
  
  Pour tout $s$ dans le centralisateur de $\psi$ vérifiant $s^2=1$, on note $\pi^A_G(\psi)(s)$ l'évaluation de $\pi^A_G(\psi)$ en $s$,
  c'est-à-dire que l'on remplace les représentations de $A(\psi)$  par leur trace en l'image de $s$ dans $A(\psi)$; 
  il s'agit ici d'évaluer en $s$ une combinaison linéaire de caractères de $A(\psi)$. 
  La première des propriétés est que $\pi^A_{G^*}(\psi)(s_\psi)$ est une représentation virtuelle stable.
  
   On note $\underline{H}_s$ la donnée endoscopique de $G$ de la forme $(s,H,\xi_H)$ telle que $\psi$ soit inclus dans l'image 
   de $\xi_H$; il est facile de vérifier que cette donnée est uniquement déterminée et qu'elle est elliptique.
  Considérons $\pi_H^A(\psi)$ associée à $\psi$  et à $H$ et  son évaluation en $s_\psi$,  $\pi_H^A(\psi)(s_\psi)$ : 
  c'est une représentation virtuelle  stable de $H$. On normalise les facteurs de transfert en fixant un modèle de Whittaker pour $G^*$ 
  et en supposant que $G$ est une forme intérieure pure de $G^*$ ou en utilisant les travaux de Kaletha.
  On demande que $e^K(G)\pi^A_G(\psi)(s_\psi s)$ soit le transfert de  $\pi^A_H(\psi)(s_\psi)$, où $e^K(G)$ est le signe de Kottwitz de $G$.
  En particulier si $s=1$, cette propriété dit que $\pi^A_G(\psi)(s_\psi)$ est le transfert de $\pi^A_{G^*}(\psi)(s_\psi)$ au signe de Kottwitz près.

Lorsque  $G$ est quasi-déployé, c'est-à-dire $G=G^*$,  les définitions précédentes ne suffisent pas  à déterminer $\pi^A_{G}(\psi)$.
 Pour les groupes classiques ({\sl cf.} Arthur  et Mok pour les groupes unitaires), il faut aussi utiliser l'endoscopie tordue   pour compléter la définition. 
 Il est alors démontré que $\pi^A_{G^*}(\psi)$ est une représentation unitaire de $G^*\times A(\psi)$ c'est-à-dire une combinaison
  linéaire à coefficients entiers positifs de représentations unitaires irréductibles de $G^*\times A(\psi)$.
  
  Il n'est pas très difficile de vérifier pour  les groupes classiques non quasi-déployés que  $\pi^A_G(\psi)$ est effectivement 
  aussi uniquement déterminée par les propriétés d'endoscopie décrites. Mais il n'est évidemment pas clair que 
  $\pi^A_G(\psi)$ soit une représentation unitaire de $G\times A(\psi)$, comme dans le cas quasi-déployé.
  L'un des buts  très modeste de cette note est de vérifier que ceci reste vrai pour les groupes classiques non nécessairement quasi-déployés.
En fait plus précisément, on le fait si $\psi$ est tempéré, c'est-à-dire trivial sur $\SL(2,\mathbb{C})$, cas connu sur les corps archimédiens
 mais pas entièrement pour les corps $p$-adiques. Dans ce cas, on complète totalement la classification, en particulier
  on montre que l'on a, comme attendu, une classification de Langlands et que la description des représentations cuspidales est analogue
   à celle du cas quasi-déployé faite par exemple en \cite{pourps}. La méthode est  une copie conforme de celle d'Arthur
    et s'appuie sur deux idées: par voie globale on montre que les coefficients sont nécessairement des entiers positifs et par voie locale
     avec l'aide du produit scalaire elliptique on montre que ces entiers sont de valeur absolue un. 
   Ces résultats dans le cas tempéré sont l'objet de la section \ref{castemperee}.

Ensuite  dans la section \ref{casunipotent}, on passe au cas des morphismes unipotents. 
Dans le cas $p$-adique ce sont les morphismes triviaux sur le groupe 
$\SL(2,\mathbb{C})$ et dans le cas archimédien ce sont ceux dont l'image du sous-groupe $\mathbb{C}^*$ de $W_F$ est centrale. 
On traite  d'abord le cas $p$-adique en utilisant l'involution d'Aubert-Schneider-Stuhler. Ce cas se ramène au cas tempéré car
 Hiraga a montré que cette involution commute à des signes explicites près à l'endoscopie. Le seul problème est donc de prendre en 
 compte ces signes, ce qui est fait en \ref{unipotentpadique}. Puis par voie globale, on traite ensuite le cas archimédien en 
 \ref{unipotentarchimedien}; on s'appuie sur \cite{pourhowe} qui ne traite pas le cas des groupes unitaires, on supposera donc que le 
 groupe est orthogonal, le cas unitaire peut se faire facilement directement comme on le vérifiera dans un article ultérieur.
On obtient en particulier un résultat de multiplicité 1 dans le cas  unipotent archimédien (théorème \ref{unipotentarchimedien})
dont on espère qu'il se généralise aux paramètres quelconques (voir \cite{MR3} pour l'énoncé de nos résultats sur les paquets d'Arthur
généraux des groupes classiques dans le cas archimédien et en particulier
cette question de la multiplicité un).

On ne va pas plus loin, c'est-à-dire que dans cet article, on ne traite pas le cas des morphismes généraux. 
Cela se fait facilement par voie locale en utilisant les descriptions locales qui sont faites pour les places $p$-adiques 
et que l'on fera dans un article en cours pour les archimédiennes. 
Mais le résultat final est bien que $\pi^A_G(\psi)$ est une représentation unitaire de $G\times A(\psi)$.

\

Les méthodes de cet article sont très fortement inspirées (voire copiées) de celles que J. Arthur a introduites pour stabiliser
 la formule des traces. Cette stabilisation introduit des objets locaux et globaux, les objets locaux sont des combinaisons linéaires
  finies de traces de représentations et on  transfère ces objets locaux. Il est donc indispensable de savoir que le transfert local
   d'une combinaison linéaire stable finie de caractères est une combinaison linéaire finie de caractères. C'est le transfert local 
   spectral. Il a été démontré par Arthur dans le cas des groupes $p$-adiques en \cite{selecta}. La référence dans le cas des places 
   archimédiennes est plus difficile à trouver car cela ne résulte pas du transfert tempéré démontré par Shelstad. Or  \cite{MW2}
    nécessite  ce résultat dans le  cadre plus général du  transfert endoscopique tordu : 
    le transfert local spectral pour les groupes $p$-adiques y  est  redémontré  (suivant \cite{selecta}) 
    en \cite{MW2} XI.2.11  (cas des représentations elliptiques) et XI.3.1 (où le cas général est ramené  au cas des représentations elliptiques), 
   et dans   le cas archimédien la démonstration est en IV.3.3. Le transfert local spectral est donc établi en toute généralité.

\

\noindent {\sl Remerciements}.  Les auteurs remercient le référé pour sa relecture attentive du texte.
C. Moeglin a bénéficié d'excellentes conditions de travail au CIRM dans le cadre de la Chaire Morlet; elle remercie le CIRM, D. Prasad et 
V. Heiermann pour cela. D. Renard bénéficie de l'ANR FERPLAY -13-BS01-0012  et remercie le CNRS et les membres 
de cette ANR pour cette aide.

\section{$A$-paquets ou plutôt $A$-représentations}
\subsection{Au sujet de la définition des A-paquets\label{rmqsurdef}} 
 On reprend sans les répéter les notations et définitions de l'introduction, en particulier, on se limite au cas des groupes orthogonaux et unitaires qui sont des formes intérieures pures de leur forme quasi-déployée. On apporte une précision nécessaire sur le fait que 
 $\pi^A_G(\psi)$ est bien défini.
 On fixe $s$ dans le centralisateur de $\psi$. On peut supposer que $s^2=1$, ce que l'on fait. 
  On a une description évidente du commutant dans la composante neutre du groupe 
 dual en utilisant les valeurs propres de $s$. Et la donnée endoscopique associée est alors elliptique mais des éléments distincts
  $s,s'$ peuvent avoir même image dans $A(\psi)$. Cela ne se produit que si $\psi$ composé avec la représentation standard du 
  $L$-groupe a de la multiplicité. On reprend sans insister la notion de bonne parité (cf. par exemple \cite{pourps}). Le seul réel problème 
  est quand des sous-représentations de $\psi$ de bonne parité interviennent avec multiplicité. Les éléments du commutant de $\psi$ de 
  carré 1 séparent ces composantes en deux sous-ensembles. Soient $s,s'$ comme ci-dessus de même image dans $A(\psi)$ alors, il 
  est facile de voir que même si les données endoscopiques ne sont pas les mêmes, elles diffèrent par le fait qu'un  nombre
   pair de copies d'une même représentation est séparé différemment par $s$ et $s'$. Le transfert de la représentation virtuelle  stable 
   correspondante est bien le même puisque le transfert commute à l'induction parabolique.

 \

  Cela permet donc de définir de façon unique $\pi^A_G(\psi)( s)$ en tout élément $s$ de $A(\psi)$ et ainsi $\pi^A_G(\psi)$ 
  est une combinaison linéaire de représentations de $G$ à coefficients dans les fonctions sur $A(\psi)$. 
  Comme $A(\psi)$ est un $2$-groupe, ceci est équivalent à dire que $\pi^A_G(\psi)$ est une combinaison linéaire à coefficients
   complexes de représentations de $G\times A(\psi)$. 
   
   \
   
   Soit $\epsilon$ un caractère de $A(\psi)$, on définit $\pi(\psi)_\epsilon$ par la formule d'inversion
   $$
   \pi(\psi)_\epsilon:=|A(\psi)|^{-1}\sum_{s\in A(\psi)} \epsilon(s) \pi^A_G(\psi)(s).
   $$
   On simplifiera un peu en \ref{hasse} (*) ci-dessous.

 \subsection{Invariant de Hasse et action de  $Z(\widehat G)^\Gamma$ \label{hasse}}
 Les formes intérieures des groupes orthogonaux (resp. unitaires) sont classifiées par les formes bilinéaires symétriques 
 (resp. hermitiennes) non dégénérées. 
 
 A une forme bilinéaire symétrique est associé un invariant de Hasse, de façon classique cet invariant n'est pas un pour un plan hyperbolique. 
 On le normalise pour remédier à cela de la façon suivante: on fixe une forme bilinéaire symétrique non dégénérée maximalement 
 isotrope sur un espace vectoriel de même dimension. L'invariant de Hasse normalisé de toute forme symétrique non dégénérée de
  même dimension est alors, par définition,
   le produit de l'invariant de Hasse classique par l'invariant de Hasse classique de la forme fixée.
   Et c'est cet invariant de Hasse normalisé que l'on appelle invariant de Hasse. Il vérifie la formule de produit local/global quand tout
    est défini sur un corps de nombres. Sur les corps $p$-adiques, les formes orthogonales sont classifiées par leur dimension, 
    leur discriminant et leur invariant de Hasse avec une exception en dimension deux. En dimension deux, on a la forme orthogonale
     déployée correspondant au discriminant un et à l'invariant de Hasse un et pour chaque discriminant non trivial, il y a deux formes
      orthogonales l'une d'invariant de Hasse un et l'autre d'invariant de Hasse moins un. Sur le corps des réels, les formes orthogonales 
      sont classifiées par leur signature ce qui est plus fin que le discriminant et l'invariant de Hasse. Même si cela n'est sans doute pas 
      indispensable, on suppose comme c'est légitime que la forme quasi-déployée $G^*$ correspond à une forme orthogonale 
      d'invariant de Hasse un.

 Pour les formes hermitiennes, la notion d'invariant de Hasse est moins classique. Notons $E$ l'extension quadratique du corps de base 
 qui définit les groupes unitaires de la situation. Pour tout élément de $E$ on définit sa partie réelle comme la demi-somme de cet élément 
 avec son conjugué. L'invariant de Hasse associé à une forme unitaire est, par définition, l'invariant de Hasse de la    forme bilinéaire symétrique 
 associée à la partie réelle de la forme hermitienne. On a donc aussi une formule de produit local/global quand on 
 part d'un corps de nombres. On suppose aussi que $G^*$ est le groupe unitaire d'une forme hermitienne d'invariant de Hasse un.
  
 Un invariant de Hasse est un signe et il définit donc un caractère de $\{\pm 1\}$ que l'on note $\epsilon_G$, c'est un peu abusif car cet invariant 
  dépend de la forme bilinéaire qui définit le groupe et  de la forme intérieure fixée représentant $G^*$.

  On fixe $\psi$ et on a donc défini $\pi^A_G(\psi)$ comme une combinaison linéaire à coefficients complexes  de représentations 
  de $G\times A(\psi)$. On remarque aussi que dans les cas qui nous occupent (groupes orthogonaux)
  $Z(\widehat G)^\Gamma=Z(\widehat G)\simeq \{\pm 1\}$, où $Z(\widehat G)^\Gamma$ désigne le sous-groupe de 
  $Z(\widehat G)$ des éléments invariants sous l'action du groupe de Galois utilisée pour définir ${}^LG$.
  
  \begin{prop} Pour tout $z\in Z(\widehat G)^\Gamma$, et pour tout $s\in A(\psi)$, on a $\pi^A_G(\psi)(sz)=
  \epsilon_G(z)\pi^A_G(\psi)(s)$.
  \end{prop}
  
 \dem   On fixe $s$ dans le centralisateur de $\psi$ tel que $s^2=1$  et on note $z$ l'élément non trivial de $Z(\widehat G)^\Gamma$. 
   A $s$ et $\psi$ correspond l'unique donnée endoscopique elliptique de $G$, $\underline{H}_s$. A $sz$ et $\psi$ correspond 
   une donnée endoscopique équivalente. Mais l'action de $z$ entraîne pour le transfert la multiplication par un signe (cf. \cite{arthurnote}, 
   fin de l'introduction, le résultat est plutôt dû à Kottwitz et Shelstad). Et on doit vérifier que ce signe est trivial exactement quand   
   l'invariant de Hasse l'est aussi. Cela vient de la normalisation des facteurs de transfert pour les formes intérieures pures de $G^*$.
  \qed
  
  \begin{cor}
  Soit $\epsilon$ un caractère de $A(\psi)$. Alors $\pi(\psi)_\epsilon$ vérifie  la formule d'inversion
   $$
   \pi(\psi)_\epsilon:=|A(\psi)/Z(\widehat{G})^\Gamma |^{-1}\sum_{s\in A(\psi)/Z(\widehat{G} )^\Gamma} \epsilon(s) \pi^A_G(\psi)(s). \eqno(*)
   $$
 \end{cor}
 Ici, il y a un abus de notation, il faut comprendre que l'on somme sur un système de représentants de  $A(\psi)/Z(\widehat{G})^\Gamma$ 
 dans $A(\psi)$ et la proposition montre que les termes de la  somme ne dépendent  pas du choix de ces représentants.

  \subsection{Conséquences de la  formule des traces\label{positiviteun}}
  
  On rappelle la stabilisation de la partie discrète spectrale de la formule des traces. Cela a été fait en \cite{arthur3} et repris 
  en \cite{MW2}. Ici la situation est globale et $G,G^*$ sont des groupes réductifs  définis sur un corps de nombres, $F$,
   et on suppose que $G^*$ est la forme quasi-déployée de $G$. On note $n^*$ la dimension de la représentation
    standard du $L$-groupe de $G$.
  
  On fixe $V$ un ensemble fini de places de $F$ qui contient les places dites ``ramifiées'', c'est-à-dire les places archimédiennes, 
  les places de petites caractéristiques et les places où il y a de la ramification. On note $G(F^V)$ le groupe des points de $G$ sur les 
  adèles de $F$ hors les places dans $V$ et $G(F_V)$ le groupe produit des points de $G$ en les complétés de $F$ en les places dans $V$. On considère l'algèbre de Hecke sphérique pour $G(F^V)$ et un caractère $c^V$ de cette algèbre 
  de Hecke. On considère $I_{disc}^G(c^V)$ la partie spectrale discrète sur laquelle les fonctions non ramifiées hors de $V$ agissent par
   le caractère $c^V$. Cette distribution vue comme distribution sur $G(F_V)$ a été stabilisée avec les références déjà données et il n'est pas
    difficile de voir que cette distribution est nulle sauf si le transfert endoscopique tordu de $c^V$ au groupe linéaire $\GL(n^*)$ est un
     caractère de l'algèbre de Hecke sphérique hors de $V$ de ce groupe correspondant à un paramètre global d'Arthur. 
     En d'autres termes, il existe un paramètre global $\psi$ d'Arthur tel qu'en toute place $v$, on obtienne par localisation un paramètre local
      $\psi_v$. On dit que $\psi$ le paramètre global est régulier s'il ne se factorise pas par un sous-groupe de Levi de $G^*$. 
      On note $\epsilon^A(\psi)$ le caractère de $A(\psi)$ associé à $\psi$ par Arthur. On fixe aussi un caractère infinitésimal à toutes les places archimédiennes. Si $G$ est un groupe spécial orthogonal pair et si $\psi$ est une somme de paramètres élémentaires tous orthogonaux pairs, on suppose qu'en une place archimédienne, le caractère infinitésimal n'est pas stable sous l'action du groupe orthogonal.
  
  Pour un  ensemble de caractères $\epsilon_v$ de $A(\psi_v)$,  $v\in V$ on considère la restriction de  $\prod \epsilon_v$ à $A(\psi)$  et on dit que l'ensemble des caractères $\{\epsilon_v\}$ est admissible si cette restriction à $A(\psi)$ est le caractère $\epsilon^A(\psi)$.

  \begin{thm} Avec les hypothèses et notations ci-dessus,
  
  $$
  \sum_{\{\epsilon_v\}}\otimes_{v\in V}\; \pi(\psi_v)_{\epsilon_v}
  $$où la somme porte sur l'ensemble des caractères admissibles,
  est une représentation unitaire de $G(F_V)$.
  \end{thm}
  Le résultat est donc que la somme qui a priori est une combinaison linéaire à coefficients complexes de représentations de $G(F_V)$
   est en fait une combinaison linéaire à coefficients des entiers positifs de représentations unitaires irréductibles de $G(F_V)$. 
   Ce théorème, comme on va le voir, est évidemment copié des résultats de \cite{book}.
  Bien sûr, cela ne montre pas la même propriété en chaque place locale ce qui est pourtant ce qui nous intéresse. 
  Pour avoir cette propriété il faudra pouvoir isoler une place et savoir que les coefficients aux autres places sont 1. 
  
  \
  
 \dem  On écrit la stabilisation de $I^G_{disc}(c^V)$. On a une somme sur toutes les données endoscopiques elliptiques non ramifiées
   hors de $V$, donc un nombre fini. De plus seules interviennent les données endoscopiques elliptiques ayant un caractère
    de l'algèbre de Hecke sphérique hors de $V$ qui se transfère en $c^V$. Cela force le $s$ de la donnée endoscopique à 
    être (à conjugaison près) dans le centralisateur de $\psi$; on suppose dès le départ, comme c'est loisible que $s^2=1$. 
    Comme on fixe $\psi$,  on suppose que  $s$ est dans le centralisateur de $\psi$.
     On a donc une somme avec des coefficients explicites des $SI^{\underline{H}_s}_{disc}(c^V)$.
  
   L'hypothèse de régularité sur $\psi$ simplifie les formules car elle  évite les termes venant des sous-groupes de Levi et les opérateurs 
   d'entrelacement associés. Sous cette hypothèse simplificatrice, ces termes stables ont été calculés par Arthur dans \cite{book} et 
   récrits par Taïbi en \cite{taibi}: les constantes explicites ne dépendent que de $^LG$ et de $\underline{H}_s$ donc sont les mêmes que pour 
   le groupe déployé et on a donc le même calcul qu'Arthur mais c'est plus simple d'aller regarder le calcul explicite fait par Taïbi où le 
   caractère infinitésimal ne joue aucun rôle. Le transfert local a été défini ci-dessus en termes précisément des $\pi(\psi)_{\epsilon_v}$ et 
   on trouve que la distribution $I^G_{disc}(c^V)$ est la formule de l'énoncé (cf. \cite{taibi} 4.0.1). Pour les groupes spéciaux orthogonaux pairs, quand le paramètre n'est pas stable sous le groupe orthogonal, on utilise l'hypothèse sur le caractère infinitésimal pour séparer aussi les représentations sous l'action du groupe orthogonal.
   D'où le résultat.
  \qed
  
  \section{Paramètres tempérés $p$-adiques}
  \label{castemperee}
  \subsection{Classification des séries discrètes des groupes orthogonaux non quasi-déployés aux places $p$-adiques}\label{castemperee}
  
  On revient à la situation locale et on s'intéresse aux morphismes de Langlands discrets qui sont des cas particuliers
   des morphismes d'Arthur. Ici le problème est le cas des groupes $p$-adiques non quasi-déployés puisque tous les autres cas sont connus.
    En particulier aux places archimédiennes, à toute série discrète (et même à toute représentation tempérée) est associée un morphisme 
    de Langlands discret et un caractère 
    du groupe des composantes du centralisateur du morphisme de Langlands  associé. Cette application est injective mais n'est pas surjective, 
    pour avoir une application surjective, il faut regarder toutes les formes intérieures pures. Et la restriction de ce caractère
     à $Z(\widehat G)^\Gamma$ correspond à l'invariant de Hasse (c'est-à-dire est trivial exactement quand l'invariant de Hasse est trivial).

  Aux places $p$-adiques quand $G$ est  quasi-déployé, on a une bijection entre les séries discrètes et les couples formés d'un paramètre de 
  Langlands discret et un caractère 
  du centralisateur de ce paramètre dont la restriction à $Z(\widehat G)^\Gamma$ est triviale (cf.\cite{book}, 
  dès l'introduction et aussi \cite{pourps} qui donne une construction explicite du morphisme de Langlands). 
  On va démontrer le même résultat dans le cas non quasi-déployé. Commençons par les énoncés:

   \begin{thm} (Arthur dans le cas quasi-déployé \cite{book} 6.5)
  Avec les notations précédentes, $\pi(\phi_v)_\epsilon$ est une représentation irréductible 
  (en particulier non nulle) et l'application $(\phi_v,\epsilon)\mapsto \pi(\phi)_\epsilon$ est bijective.
  \end{thm}

  \begin{prop}
  Soit $G_v$ un groupe spécial orthogonal ou unitaire non nécessairement quasi-déployé et soit $\phi_v$ un morphisme
   de Langlands borné et discret. Alors la représentation de $G_v\times A(\phi_v)$ associée est de la forme:
  $$
  \epsilon^K(G_v) \sum_{\epsilon} \pi(\phi_v)_\epsilon \otimes \epsilon,
  $$
  où $\epsilon$ parcourt l'ensemble des caractères de $A(\phi_v)$ de restriction fixée à $Z(\widehat G)^\Gamma$ par l'invariant 
  de Hasse et où $\pi(\phi_v)_\epsilon$ est une série discrète  irréductible. 
  \end{prop}

La preuve est la même que dans le cas quasi-déployé et donc grandement due à Arthur. La méthode est la suivante: 
on commence par démontrer par voie globale que les $\pi(\phi_v)_\epsilon$ sont des sommes à coefficients entiers positifs 
de représentations. Ensuite on utilise les relations d'orthogonalité.
 Détaillons cela en plusieurs sous-sections. Pour simplifier les notations on se limite au cas des groupes orthogonaux. 
  Le cas des groupes unitaires est déjà connu (cf. \cite{kalethaandco}).

\subsection{Positivité\label{positivite}}
Cette partie nécessite de  globaliser la situation.
On commence par fixer un paramètre $\phi_v$ borné (tempéré) pour le groupe local $G_v$. 
Donc on ne suppose pas que le paramètre local est discret. Montrons alors le lemme:
\begin{lem}
Les représentations $\pi(\phi_v)_\epsilon$ sont des combinaisons linéaires à coefficients positifs de représentations irréductibles.
\end{lem}

  On décompose $\phi_v$ en somme de représentations irréductibles de $W_{F_v}\times \SL(2,\mathbb{C})$. 
  Chaque composante irréductible est une représentation symplectique si $G_v$ est un groupe orthogonal impair et une 
  représentation orthogonale sinon. Notons $m$ la dimension d'une telle composante irréductible. On va globaliser chacune 
  de ces composantes irréductibles en construisant une représentation cuspidale unitaire de $\GL(m,\mathbb{A}_F)$ ayant en la place $v$ 
  la série discrète définie par cette composante irréductible et on veut aussi des conditions aux places archimédiennes. 
  
  On construit un corps global $F$ avec une place $v$ où $F$ se localise en $F_v$ et un nombre 
  de places réelles grand par rapport au nombre de sous-représentations irréductibles incluses dans $\phi$. 
  
  En chaque place réelle on  suppose que le caractère infinitésimal est formé de demi-entiers si $G$ est un groupe orthogonal de dimension impair
   et  d'entiers sinon et on suppose que ces nombres sont tous distincts et ils seront distincts pour chaque construction: pour faire cela,
    il faut commencer par regarder un groupe classique quasi-déployé de même nature que $G$ sur un espace vectoriel de dimension $m$. 
    Le paramètre local définit un paquet de séries discrètes pour ce groupe, avec une toute petite difficulté dans le cas où $m$ est un entier impair, 
    il faut supposer  que le morphisme composé avec le déterminant soit trivial. 
    On s'y ramène en tordant par un caractère de $W_{F_v}$. On prend une série discrète dans ce paquet local. Aux places à l'infini
     on construit une série discrète pour le groupe réel avec un caractère infinitésimal que l'on choisit comme on veut et on sait 
     que l'on peut globaliser ces constructions, c'est-à-dire construire une représentation cuspidale du groupe sur les
      adèles de $F$ ayant en la place finie et en les places archimédiennes les constructions faites. 
       La théorie d'Arthur associe à cette représentation un paramètre qui est symplectique en toute place si la composante irréductible
        fixée était symplectique et qui est orthogonal sinon et un transfert tordu à $\GL(m,\mathbb{A}_F)$. Si on a tordu par un caractère 
        quadratique le paramètre local on retord globalement. Cela ne change pas le caractère infinitésimal aux places archimédiennes.

  En faisant cela pour toutes les composantes irréductibles de $\phi_v$, on obtient bien un paramètre global $\phi$ 
  tel que $A(\phi)$ s'envoie surjectivement sur $A(\phi_v)$. Ce paramètre $\phi$ est globalement discret, il a un 
  caractère infinitésimal régulier c'est-à-dire que les représentations cuspidales  des $\GL(m,\mathbb{A})$ construites 
  sont toutes distinctes. En toute place $w$, le déterminant de $\phi_w$ 
  (la localisation de $\phi$ en la place $w$) fixe le discriminant de la forme orthogonale que l'on va considérer. 
  Ce discriminant est un si $G_v$ est un groupe orthogonal impair. 
  On n'a pas encore fixé l'invariant de Hasse de la forme orthogonale non quasi-déployée qui va nous intéresser, 
  elle dépend des choix que l'on va faire maintenant.
    
   On fixe aussi un caractère $\epsilon_v$ de $A(\phi_v)$, où $v$ est la place fixée $p$-adique qui nous intéresse. 
   On a construit les objets aux places archimédiennes de sorte que l'application de $A(\phi)$ dans l'analogue en
    ces places soit injectif. On peut ainsi trouver en ces places, $w$, un
    caractère $\epsilon_w$ tel que l'ensemble qu'ils forment avec $\epsilon_v$ soit admissible c'est-à-dire ici trivial en
     restriction à $A(\phi)$ (le caractère d'Arthur est trivial pour les morphismes globaux triviaux sur le $\SL(2,\mathbb{C})$ dit de
      Lefschetz, mais peu importe ici).  En toute place $w$ comme précédemment, on considère une forme orthogonale dont l'invariant de 
     Hasse est la restriction à $Z(\widehat G)^\Gamma$ du caractère $\epsilon_w$. En la place $v$ et en les places archimédiennes, 
     on a donc le discriminant et l'invariant de Hasse d'une forme orthogonale, le produit des invariants de Hasse vaut 1.   
     Cela ne fixe pas encore la forme orthogonale aux places archimédiennes, et on en fixe une de sorte que le caractère
      $\epsilon_w$ corresponde bien à une représentation tempérée du groupe orthogonal associé. On complète aux autres
       places finies en prenant des formes orthogonales déployées avec le discriminant déjà  fixé, ce qui est loisible
        puisque le produit des invariants de Hasse vaut bien 1.
  
  On utilise maintenant \ref{positiviteun} où $\psi$ est $\phi$ prolongé trivialement à $\SL(2,\mathbb{C})$ (on garde la notation $\phi$). 
  On voit que dans la somme on a un terme 
  $\pi(\phi_v)_{\epsilon_v}\otimes \Pi^v$, où $\Pi^v$ est un produit en toute place dans $V$ sauf $v$ de représentations irréductibles,
   plus   d'autres termes de la forme $ \pi(\phi_v)_{\epsilon'_v} \otimes \Pi^{'v}$ mais où les $\Pi^{'v}$ sont eux aussi irréductibles
    mais non isomorphes à $\Pi^v$, c'est le point important et ceci est dû au fait que pour toutes les places autres que $v$ on sait que
     le caractère du centralisateur détermine au plus une représentation irréductible. Ainsi il ne peut y avoir de simplification entre  
  $\pi_{\epsilon_v}\otimes \Pi^v$ et les autres termes.
   Cela force $\pi_{\epsilon_v}$ à être une somme de représentations unitaires avec des coefficients entiers positifs. Et cela termine la preuve.
  \qed
  
  \subsection{Multiplicité un}
  On se replace dans le cas local $p$-adique. On fixe encore un morphisme de Langlands, que l'on suppose ici discret, $\phi$ de $W'_F$ dans $^LG$. 
  A ce morphisme correspond une somme  stable de séries discrètes pour $G^*$ et on a montré que son transfert à $G$ est, au signe 
  de Kottwitz près, une somme à coefficients entiers positifs de représentations unitaires, $\pi^A_G(\phi)(1)$ où $1$ est l'unité de $A(\phi)$. 
  Comme le transfert commute au module de Jacquet, le critère de Casselman pour caractériser les séries discrètes s'applique
   à tous les termes de $\pi^A_G(\phi)(1)$   et cette représentation est donc une combinaison linéaire à coefficients entiers positifs 
   de séries discrètes. On rappelle qu'à tout caractère $\epsilon$ du centralisateur de $\phi$ dont la restriction à $Z(\widehat G)^\Gamma$
    est l'invariant de Hasse, on a associé une représentation $\pi(\phi)_\epsilon$ qui est une combinaison linéaire à coefficients 
    positifs de représentations unitaires et la somme de ces représentations est $\pi^A_G(\phi)(1)$. Chacune de ces
     représentations $\pi(\phi)_\epsilon$ est donc une combinaison linéaire de séries discrètes.
  
  On va maintenant utiliser  les relations d'orthogonalité pour le produit scalaire elliptique pour prouver l'irréductibilité des 
  représentations $\pi(\phi)_\epsilon$ et le fait que ces représentations sont inéquivalentes.
  
   On reprend \cite{book} 6.5 bien que le cas qui y est traité soit  plus compliqué que ce que l'on fait ici car il utilise aussi l'endoscopie
    tordue. On fait aussi remarquer au lecteur que le fait que l'on ait démontré directement que les représentations qui nous intéressent
     sont des séries discrètes, facilite  l'usage du produit scalaire, car la norme elliptique d'une série discrète est un, ce qui n'est pas vrai
      pour une représentation elliptique générale. Pour le produit scalaire elliptique,  on renvoie à \cite{LMW} 4.6 qui reprend en détail ces questions.
  
  On commence par inverser les formules donnant $\pi(\phi)_\epsilon$. On remarque que pour $s$ fixé dans $A(\phi)$, 
  en notant $s_0$ un élément de $Z(\widehat G)^\Gamma$, on a montré l'égalité
  $$
  \epsilon(s) \, \mathrm{transfert}\left(\pi_{st}^{\underline{s,H,\xi}}(\phi_s)\right)=
  \epsilon(ss_0) \, \mathrm{transfert}\left(\pi_{st}^{\underline{ss_0,H,\xi}}(\phi_{ss_0})\right),$$ 
  où l'on a noté $\phi_s$ et $\phi_{ss_0}$ les factorisations de $\phi$ par les $L$-groupes venant de la donnée endoscopique. 
  
  On a donc la formule d'inversion
  $$
  e^K(G)\pi(\phi)_\epsilon=n(\phi)^{-1}\sum_{s\in A(\phi)/Z(\widehat{G})}\epsilon(s)\, \mathrm{transfert}\left(\pi_{st}^{\underline{s,H,\xi}}(\phi_s)\right),
  $$où $n(\phi)$ est le cardinal du groupe $A(\phi)/Z(\widehat{G})$.

  Soit $\epsilon,\epsilon'$ deux caractères de $A(\phi)$ de restrictions à $Z(\widehat G)^\Gamma$ égales à l'invariant de Hasse. 
  On calcule le produit scalaire elliptique 
  $$
  \langle \pi(\phi)_\epsilon,\pi(\phi)_{\epsilon'}\rangle_{ell}. \eqno(1)
  $$
  Les distributions obtenues par transfert pour deux données endoscopiques inéquivalentes sont orthogonales et (1) vaut donc
  $$n(\phi)^{-2}
  \sum_{s\in A(\phi)/Z(\widehat{G})}\epsilon(s)\epsilon'(s) \langle
  \mathrm{transfert}\left(\pi_{st}^{\underline{s,H,\xi}}(\phi_s)\right),\mathrm{transfert}\left(\pi_{st}^{\underline{s,H,\xi}}(\phi_s)\right)
  \rangle_{ell}.\eqno(2)
  $$On rappelle que le transfert commute au produit scalaire elliptique. Ainsi dans (2) ci-dessus, on peut remplacer $\langle
  \mathrm{transfert}\left(\pi_{st}^{\underline{s,H,\xi}}(\phi_s)\right),\mathrm{transfert}\left(\pi_{st}^{\underline{s,H,\xi}}(\phi_s)\right)
  \rangle_{ell}$ par son analogue en enlevant le mot transfert. Le calcul se fait alors dans le cas quasi-déployé et ce produit scalaire vaut exactement le nombre de séries discrètes dans le transfert pour le cas quasi-déployé, c'est-à-dire encore le nombre de séries discrètes dans le paquet associé à $\phi$ dans le cas quasi-déployé, c'est-à-dire encore $n(\phi)$.
  Ainsi (2) est nul si $\epsilon\neq \epsilon'$ et si $\epsilon=\epsilon'$ (2) vaut un.
   On a donc les conséquences suivantes: les $\pi(\phi)_\epsilon$ sont des représentations disjointes quand $\epsilon$ varie. 
   Et leur norme est aussi la norme est un. Il n'est pas non plus difficile de voir que $\pi(\phi)_\epsilon$ est orthogonal
    à $\pi(\phi')_\eta$ si $\phi'$ est non conjugué de $\phi$ pour les mêmes raisons. 
  
  Il reste à montrer que toute  série discrète de $G$ est bien l'une des représentations $\pi(\phi)_\epsilon$ que nous venons d'étudier.
   Pour cela on utilise encore la stabilisation de l'espace des fonctions $I_{cusp}(G)$ identifiée à l'espace des pseudo-coefficients des séries 
   discrètes de $G$. Soit $\pi$ une série discrète pour $G$. Un argument de Waldspurger déjà utilisé en \cite{pourps} 2.4, montre que la projection
    orthogonale sur la partie stable de $I_{cusp}(G)$ d'un pseudo-coefficient de $\pi$ est non nulle. Ainsi il existe un paquet stable de
     représentations elliptiques de $G^*$ dont le transfert à $G$ n'est pas orthogonal à $\pi$ pour le produit scalaire elliptique. 
     Les paquets stables de représentations elliptiques pour $G^*$ sont des sommes de séries discrètes et il existe donc $\phi$ 
     un paramètre de Langlands comme ceux étudiés ci-dessus, tel que $\pi$ ne soit pas orthogonal au transfert du paquet stable
      de $G^*$ associé à $\phi$. Or ce transfert n'est autre que $\sum_\epsilon \pi(\phi)_\epsilon$. Ainsi $\pi$ est l'une des représentations
       $\pi(\phi)_\epsilon$. Cela termine la preuve du théorème et de la proposition \ref{castemperee}.

  \subsection{Au sujet des modules de Jacquet des séries discrètes\label{Jac}}
  Ce paragraphe ne nous servira pas pour le but principal de ce papier mais est mis ici par souci de complétion.

  Dans \cite{europe} et \cite{ams}, on avait montré, avec M. Tadic, que les modules de Jacquet des séries discrètes
   des groupes classiques $p$-adiques donnaient des renseignements précis sur les paramètres de Langlands. 
   Et en particulier que ces propriétés de modules de Jacquet ramenaient la paramétrisation de ces séries discrètes 
   au cas cuspidal et on décrivait les paramètres des représentations cuspidales sans montrer dans ces références 
   que tous ces paramètres potentiels correspondaient bien à une représentation cuspidale. Cela a été fait ultérieurement
    (cf. \cite{pourps}) en utilisant les méthodes de \cite{selecta} (qui sont celles utilisées ici aussi). 
    Comme \cite{pourps} utilisait le fait que les groupes soient quasi-déployés, on va ici redonner l'argument pour les groupes
     non quasi-déployés pour que ceux que cela intéresse puissent avoir une référence complète même si il n'y a rien de nouveau.
  
  On fixe $\phi$ un paramètre de Langlands discret pour $G$ et on reprend la notation $Jord(\phi)$ pour les sous-représentations
   irréductibles de $W_F\times \SL(2,\mathbb{C})$ incluses dans $\phi$. On identifie les éléments de 
   $Jord(\phi)$ à des couples $(\rho,a)$ où $\rho$ est la représentation cuspidale unitaire d'un groupe $GL(d_\rho,F)$
    (ce qui définit $d_\rho$) qui correspond à la représentation irréductible de $W_F$ et $a$ la dimension de la 
    représentation irréductible de $\SL(2,\mathbb{C})$.

  Pour $\rho$ comme ci-dessus, $x\in \mathbb{C}$ et $\pi$ une représentation irréductible de $G$, 
  on reprend  la notation $Jac_{\rho||^x}\pi$ pour noter la composante $\rho||^x$ isotypique pour l'action de $GL(d_\rho,F)$ 
  dans le module de Jacquet de $\pi$ pour le parabolique (standard) maximal de $G$ ayant un sous-groupe de Levi contenant 
  un facteur $GL(d_\rho,F)$, s'il en existe; sinon $Jac_{\rho||^x}\pi$ est par définition $0$. 
  Ainsi $Jac_{\rho||^x}(\pi)$ est une représentation du groupe $\GL(d_\rho,F)\times G'$ où $G'$ 
  est un groupe de même type que $G$ mais de rang plus petit dont tous les sous-quotients irréductibles 
  sont de la forme $\rho||^x\times \sigma$ où $\sigma$ est une représentation cuspidale d'un sous-groupe
   de Levi de $G'$. Quand la représentation $Jac_{\rho||^x}\pi$ est irréductible, on la voit simplement comme une représentation de $G'$.
  
  Soit $\phi$ comme ci-dessus et soit $\epsilon$ un caractère de $A(\phi)$ dont la restriction à $Z(\widehat G)^\Gamma$ correspond à l'invariant de 
  Hasse de la forme orthogonale que $G$ préserve. Et on note $\pi_\epsilon$ la série discrète correspondant à $\phi$ et 
  $\epsilon$, comme $\phi$ est fixé, on ne le met pas dans la notation. On identifie $\epsilon$ à une application de 
  $Jord(\phi)$ dans $\{\pm 1\}$.
  
  \begin{thm} Soit $\rho$ et $x$ comme ci-dessus. Alors
  $Jac_{\rho||^x}\pi_\epsilon=0$ sauf éventuellement si $x>0$ et il existe $a\in \mathbb{N}$ tel que $(\rho,a)\in Jord(\phi)$ 
  avec $x=(a-1)/2$ (d'où nécessairement $a\geq 2$).
  
  Si ces conditions sont remplies alors on a effectivement 
 $Jac_{\rho||^x}\pi_\epsilon\neq 0$ exactement quand l'une des conditions ci-dessous est satisfaite: 
 
 $a>2$ et $(\rho,a-2)\notin Jord(\phi)$ 
 
 ou  $(\rho,a-2)\in Jord(\phi)$ et $\epsilon((\rho,a))=\epsilon((\rho,a-2))$
 
 ou $a=2$ et $\epsilon((\rho,a))=+1$.
 \end{thm}
  
 \dem  On écrit 
  $$
  \pi_\epsilon=(|A(\phi)/Z(\widehat{G})|)^{-1}\sum_{s\in A(\phi)/Z(\widehat{G})}\epsilon(s)\,  \mathrm{transfert}(\pi( \phi_s)(1)),
  $$
  où on a noté  $\pi( \phi_s)(1)$ la distribution stable pour la donnée endoscopique elliptique associée à 
  $s$ et à $\phi$. On utilise le fait que le transfert commute au module de Jacquet pour écrire avec les notations de l'énoncé:
  $$Jac_{\rho||^x}\pi_\epsilon=(|A(\phi)/Z(\widehat{G})|)^{-1}\sum_{s\in A(\phi)/Z(\widehat{G})}\epsilon(s)\, 
   \mathrm{transfert}(Jac_{\rho||^x}\pi( \phi_s)(1)). \eqno(1)
  $$
  Les représentations $Jac_{\rho||^x}\pi( \phi_s)(1)$ sont connues puisque l'on est dans le cas
   quasi-déployé. Toutes ces représentations sont nulles sauf éventuellement si $x=(a-1)/2$ avec $a>1$ 
   tel que $(\rho,a)\in Jord(\phi)$. Et si ces conditions sont satisfaites, ces représentations sont 
   stables,  associées au morphisme $\phi_-$ où $\phi_-$ se déduit de $\phi$ en remplaçant $(\rho,a)$ par $(\rho,a-2)$ dans la 
   décomposition de $\phi$ en sous-représentations irréductibles. En particulier si $(\rho,a-2)\notin Jord(\phi)$ et $a>2$,  $\phi_-$
    est un paramètre de séries discrètes pour le groupe qui se déduit de $G^*$ en diminuant le rang par $d_\rho$ et $A(\phi_-)$
     est naturellement isomorphe à $A(\phi)$. Dans ce cas, par définition, on a alors $Jac_{\rho||^x}\pi_\epsilon= \pi(\phi_-)_\epsilon$ si $x=(a-1)/2$.  
  
  Il reste donc à voir le cas où soit $a=2$ soit  $(\rho,a-2)\in Jord(\phi)$ avec $x=(a-1)/2$. Si $a=2$, $\phi_-$ 
  s'obtient en enlevant $(\rho,a)$. Dans les deux cas, 
  l'application naturelle de $A(\phi)$ dans $A(\phi_-)$ est surjective mais non injective. On note $T$ le noyau, 
  il a deux éléments. On remarque que dans (1), on a une sous-somme sur les éléments de $T$ et que (1) est 
  non nul, seulement si le caractère $\epsilon$ est trivial sur $T$, c'est-à-dire si $\epsilon(\rho,a)=\epsilon(\rho,a-2)$ 
  avec la convention que $\epsilon(\rho,a-2)=+1$ si $a=2$. Ce sont les conditions de l'énoncé. Si ces conditions 
  sont satisfaites, alors le terme de gauche est par définition $\pi(\phi_-)_\epsilon$ où $\epsilon$ est maintenant vu 
  comme un caractère de $A(\phi_-)$. Cela termine la preuve.\qed

  \subsection{Classification des représentations cuspidales}
  \begin{cor}
  Les représentations cuspidales de $G$ sont classifiées exactement comme dans le cas quasi-déployé: il n'y  en a pas si 
  le morphisme $\phi$ correspondant est tel que $Jord(\phi)$ a des trous, c'est-à-dire s'il existe $(\rho,a)\in Jord(\phi)$ avec $a>2$ 
  et $(\rho,a-2)\notin Jord(\phi)$. Et si le morphisme $\phi$ n'a pas de trous, les représentations cuspidales
   correspondent aux caractères de $A(\phi)$ de restriction à $Z(\widehat{G})$
   l'invariant de Hasse et qui sont alternés c'est-à-dire $\epsilon(\rho,a)\neq \epsilon(\rho,a-2)$ si $a\geq 2$ en posant 
   $\epsilon(\rho,0)=+1$ si nécessaire.
  \end{cor} 
  
    \subsection{Classification des représentations tempérées}
  
  \begin{prop} Soit $\phi$ un morphisme tempéré de $W'_F$ dans $^LG$. 
  Pour tout caractère $\epsilon$ de $A(\phi)$ de restriction  à $Z(\widehat{G})$  l'invariant de Hasse, la représentation 
  $\pi(\phi)_\epsilon$ est irréductible et tempérée. L'application qui à $(\phi,\epsilon)$ associé $\pi(\phi)_\epsilon$ est une bijection.
  \end{prop}
  
 \dem  Il est clair par leur définition même en tant que transfert que les représentations $\pi(\phi)_\epsilon$ sont
   tempérées (cf. \cite{MW2} XI.2.11).  On considère la combinaison linéaire: 
  $$
  \sum_\epsilon \pi(\phi)_\epsilon.\eqno(1)
  $$
  C'est, par définition, le transfert de la représentation tempérée stable associée à $\phi$ pour $G^*$. 
  On connaît ces représentations stables pour $G^*$: ce sont des induites complètes d'un paquet stable 
  de séries discrètes pour un sous-groupe de Levi, $M^*$ de $G^*$. Ainsi il existe un sous-groupe de Levi 
  $M^*$ de $G^*$ tel que $\phi$ se factorise en un paramètre de séries discrètes pour $M^*$. Et la représentation (1) est le transfert
   de $G^*$ à $G$ de l'induite complète du paquet stable de séries discrètes de $M^*$. Ce transfert est nul si $M^*$ ne se transfère
    pas à $G$. Or d'après \ref{positivite} la somme des représentations $\pi(\phi)_\epsilon$ voit toutes les composantes de chacune 
    de ces représentations (les coefficients sont positifs et s'ajoutent). Ainsi cette somme ne peut être nulle et il faut donc que $M^*$
     se transfère. D'où $M$. D'autre part, on vient de montrer aussi que chaque composante irréductible d'une des représentations 
     $\pi(\phi)_\epsilon$ est  sous-module d'une induite à partir d'un sous-groupe parabolique de Levi $M$ et d'une série discrète de
      $M$. Et on rappelle que le support discret d'une représentation tempérée, $\pi$, est bien défini, c'est-à-dire qu'à conjugaison 
      près la donnée d'un couple formé d'un sous-groupe de Levi $M$ et d'une série discrète $\sigma$ de $M$ tel que $\pi$ soit
       sous-module de l'induite de $\sigma$ après choix d'un sous-groupe parabolique de Levi $M$ est uniquement déterminé par $\pi$.
    
    Cela montre donc deux points importants: si $\phi$ n'est pas conjugué d'un morphisme tempéré $\phi'$ alors pour tout
     $\epsilon$ et $\epsilon'$, des caractères des groupes des centralisateurs de $\phi$ et $\phi'$ respectivement, 
     $\pi(\phi)_\epsilon$ n'a aucune composante en commun avec  $\pi(\phi')_{\epsilon'}$. Et d'autre part, pour toute
      représentation tempérée irréductible, il existe $\phi$ et $\epsilon$ tel que $\pi$ soit une composante irréductible 
      intervenant dans la combinaison linéaire définissant $\pi(\phi)_\epsilon$.
    
    Il reste donc à montrer que les $\pi(\phi)_\epsilon$ sont des représentations irréductibles non isomorphes entre elles quand $\phi$ est fixé.
 Il faut donc démontrer que les représentations $\pi(\phi)_\epsilon$ sont non nulles et que les induites de séries discrètes n'ont pas de multiplicité.
 
 On commence par un exemple. On considère le cas d'un sous-groupe parabolique maximal de $G$ de la forme $\GL(d,F)\times G'$
  avec $G'$ un groupe de même type que $G$ et on fixe une série discrète de ce sous-groupe de Levi. On remarque ici que $G'$ 
  n'est pas le groupe trivial si $G$ n'est pas quasi-déployé. Alors l'induite est soit irréductible soit de longueur deux, un calcul de 
  module de Jacquet prouve  cela (cf. \cite{ams}). Le R-groupe est donc soit trivial soit isomorphe à $\mathbb{Z}/2\mathbb{Z}$. 
  La décomposition de l'induite est, par la théorie générale, gouvernée par les représentations irréductibles d'une extension
   centrale par $\mathbb{C}^*$, du R-groupe ayant un caractère central fixé. Or de telles extensions centrales pour $\mathbb{Z}/2\mathbb{Z}$ sont non seulement abéliennes (ce qui est vrai pour toute extension centrale d'un groupe cyclique) mais même scindée: en effet soit $$1\rightarrow \mathbb{C}^* \rightarrow C \rightarrow \mathbb{Z}/2\mathbb{Z} \rightarrow 1$$ une telle extension. On prend $z$ un élément de $\mathbb{Z}/2\mathbb{Z}$ et $\overline{z}$ un relèvement dans $C$. Quitte à multiplier par un élément de $\mathbb{C}$, on s'arrange pour que $\overline{z}^2=1$. Et dans ces conditions $\overline{z}$ est égal à son inverse. Ainsi $C$ qui est engendré par $\overline{z}$ et son centre est abélien et est même une extension scindée.
    L'induite est donc sans multiplicité. Et elle est  réductible si et seulement si il existe une représentation elliptique
     combinaison linéaire des sous-représentations incluses dans l'induite. Comme les représentations elliptiques s'obtiennent 
     par transfert endoscopique (cf. \cite{MW2} XI.4) on a une description précise de ce cas.
 
 Pour avoir le cas général, on se heurte à la difficulté que pour les places $p$-adiques on ne sait pas que la décomposition
  des induites de séries discrètes est gouvernée par les représentations du R-groupe lui-même, il semble même que ce ne 
  soit pas vrai. On prend donc une autre méthode, celle utilisée pour ramener le cas d'un paramètre d'Arthur général au cas
   d'un paramètre d'Arthur ``discret'' et ici cela veut vraiment dire discret. Et on démontrera que ce sont bien les représentations
    du R-groupe lui-même qui gouvernent la décomposition des induites.
 
 On décompose $\phi$ en sous-représentations irréductibles de $W_F\times \SL(2,\mathbb{C})$ et l'on écrit $\phi=\oplus_{i\in [1,\ell]}
  \rho_i\boxtimes R[a_i]$ où les $a_i$ sont ordonnées tels que $a_1\geq \cdots\geq  a_\ell$. On fixe des entiers $T_1 >> \cdots >>T_\ell$
   et on considère le morphisme $\phi_+$ pour un groupe de même type que $G$ mais de rang $\sum_i T_i$ plus grand dont la 
   décomposition en sous-représentations irréductibles est:
 $$
 \oplus_{i\in [1,\ell]} \rho_i\boxtimes R[a_i+2T_i].$$
 On a vérifié par des arguments assez élémentaires de modules de Jacquet (cf. \cite{pourshahidi} 3.2) que pour tout
  $\epsilon_+$ caractère de $A(\phi_+)$ la représentation
 $$
 \circ_{i\in[\ell,1]; k\in [T_i,1]}Jac_{\rho_i||^{(a_i-1)/2+k}}\; (\pi(\phi_+)_{\epsilon_+}) \eqno(2)$$
 est soit nulle soit irréductible et que si elle est non nulle elle n'est pas isomorphe à une représentation du même type 
 avec un autre caractère $\epsilon_+$. 
 On a une application naturelle de $A(\phi_+)$ dans $A(\phi)$ et on vérifie comme on l'a fait dans la preuve du théorème
  de \ref{Jac} que (2) est nul si $\epsilon_+$ ne se factorise pas en un caractère de $A(\phi)$. On note alors $\epsilon$
   cette factorisation et le fait que l'endoscopie commute aux modules de Jacquet force alors (2) à être isomorphe à $\pi(\phi)_\epsilon$.
 
 Ainsi les représentations $\pi(\phi)_\epsilon$ sont soit irréductibles soit nulles et non isomorphes entre elles. 
 On vient donc de démontrer que les induites de séries discrètes se décomposent sans multiplicité et cela entraîne
  donc que les représentations projectives du R-groupe qui gouvernent cette décomposition se relèvent en des représentations 
  du R-groupe. On fait alors remarquer au lecteur que la preuve ci-dessus, qui s'applique aussi au cas quasi-déployé est
   légèrement différente de celle de \cite{book}. En {\sl loc. cite}, Arthur calcule vraiment les opérateurs d'entrelacement
    et obtient donc plus directement la remarque ci-dessous.
 
On reprend les notations précédentes pour $\phi$ un paramètre de Langlands borné (tempéré) pour $G$.
 On note $M^*$ le sous-groupe de Levi de $G^*$ pour lequel $\phi$ définit un paquet de séries discrètes. 
 Si $M^*$ ne se transfère pas à $G$, il n'y a pas de représentations tempérées pour $G$ associées à $\phi$ 
 comme on l'a vu et sinon on note $M$ le transfert de $M^*$ et on fixe $\sigma$ une série discrète de $M$ 
 associée à $\phi$. On fixe un sous-groupe parabolique de $G$ de sous-groupe de Levi $M$ et on considère
  l'induite de $\sigma$. On note $\ell$ le nombre de sous-représentations irréductibles de $\phi$ y intervenant avec une multiplicité paire.

 \begin{rmq} Le R-groupe de l'induite de $\sigma$ est isomorphe à $(\mathbb{Z}/2\mathbb{Z})^\ell$ et les 
 caractères de ce groupe sont en bijection avec les sous-quotients irréductibles de l'induite conformément à la théorie du R-groupe. 
 \end{rmq}
 Pour avoir ce résultat unifié, il faut considérer le groupe orthogonal pair au lieu du groupe spécial orthogonal.
 
  Le fait que le R-groupe soit isomorphe à $(\mathbb{Z}/2\mathbb{Z})^\ell$ résulte du cas où 
$M$ est maximal, cas que l'on a vu ci-dessus en exemple, et des résultats de Goldberg (\cite{goldberg}, 6.5 qui 
expliquent la subtilité pour $\SO(2n)$). Comme on a montré a priori le fait que l'induite se décompose avec multiplicité
 un, on sait que l'extension centrale du R-groupe qui gouverne la décomposition de cette induite est abélienne. 
 Comme on l'a noté plus haut,  les extensions centrales  par $\mathbb{C}^*$d'un produit de groupes $
  \mathbb{Z}/2\mathbb{Z}$ sont scindées et abéliennes. D'où la remarque.

   Il nous restait à montrer que chaque $\pi(\phi)_\epsilon$ est non nul.  Le plus rapide est d'utiliser la remarque pour compter le nombre 
   de sous-représentations dans l'induite. En fixant $\sigma$ comme dans la remarque, on fixe la restriction de $\epsilon$ à $A_M(\phi)$,
    l'analogue de $A(\phi)$ pour $M$. On note $\epsilon_M$ cette restriction. Le nombre de caractères du R-groupe est alors exactement
     le nombre d'extensions de $\epsilon_M$ en un caractère de $A(\phi)$. Et toutes les extensions ont nécessairement comme restriction 
     à $Z(\widehat G)^\Gamma$ la  valeur de $\epsilon_M$ sur  $Z(\widehat M)^\Gamma$, donc l'invariant de Hasse comme requis.
      Cela termine la preuve. \qed

  \section{Le cas unipotent}\label{casunipotent}
  \subsection{Le cas $p$-adique\label{unipotentpadique}}
  On fixe encore un corps local $p$-adique.
  On dit qu'un morphisme $\psi$ est unipotent si sa restriction au groupe de Weil-Deligne du corps local $F$ est trivial sur
   $\SL(2,\mathbb{C})$. C'est une définition commode mais qui ne correspond pas à  l'analogue de celle donnée par
    Barbasch-Vogan dans le cas archimédien. Il vaudrait sans doute mieux dire que le paramètre est anti-tempéré car 
     les représentations associées à un tel morphisme sont les duales au sens d'Aubert-Schneider-Stuhler des représentations 
     tempérées comme on va le revérifier.
  
Toutefois la dualité ne commute pas tout à fait au transfert comme démontré par Hiraga (\cite{hiraga}) et 
 introduit une permutation dans les caractères qui paramètrent les représentations. Donc avant de pouvoir 
 énoncer le théorème, il faut introduire un caractère qui va tordre la situation. On note $D$ l'involution d'Aubert-Schneider-Stuhler
  et pour tout groupe $G$, on note $D_G:=(-1)^{a_G}D$ l'involution $D$ multipliée par un signe où $a_G$ est la dimension d'un tore 
  déployé maximal de $G$. Ainsi grâce aux résultats d'Hiraga déjà cités, cette involution commute à l'endoscopie.

On fixe $\phi$ un paramètre tempéré pour $G^*$. On note $1$ le caractère trivial de $A(\phi)$ et on a donc la représentation
 tempérée $\pi(\phi)_1$.  On note $t(\phi)$ le signe tel que $t(\phi)D_{G^*}(\pi(\phi)_1)$ soit une représentation irréductible donc avec un
  signe +. Ce signe est connu.

Soit $s\in A(\phi)$, on pose $\epsilon_\phi(s):=t(\phi)t(\phi_s)$ où $t(\phi_s)$ est l'analogue de $t(\phi)$ pour la donnée endoscopique 
elliptique associée à $s$ et $\phi$ et pour le paramètre $\phi_s$ de cette donnée endoscopique qui est la factorisation de $\phi$.

\begin{lem} L'application $s\mapsto \epsilon_\phi(s)$ est un caractère de $A(\phi)$ trivial sur  $Z(\widehat{G})^\Gamma$ et sur
 l'image par $\phi$ du centre de $\SL(2,\mathbb{C})$.
\end{lem}
\dem Pour démontrer ce lemme, on se ramène au cas où $\phi$ est trivial sur $\SL(2,\mathbb{C})$. On a la propriété d'invariance
 suivante pour $t(\phi)$. Soit $a$ un entier tel que la restriction de $\phi$ à $\SL(2,\mathbb{C})$ contienne la représentation 
 irréductible de dimension finie de dimension $a$. On note $\phi_-$ la représentation de $W_F\times \SL(2,\mathbb{C})$ 
 qui se déduit de $\phi$ en remplaçant toutes les représentations irréductibles de $\SL(2,\mathbb{C})$ de dimension $a$
  par des représentations irréductibles de dimension $a-2$. En termes clairs si $\phi=\oplus_{(\rho,b)}\rho\otimes R[b]$ est
   la décomposition en sous-représentations irréductibles de $\phi$, celle de $\phi_-$ se déduit de celle de $\phi$ en remplaçant les couples $(\rho,b)$
    tel que $b=a$ par $(\rho,a-2)$. On vérifie que $t(\phi)=t(\phi_-)$: la raison est que $\pi(\phi)_1$ est un sous-quotient irréductible
     de l'induite $\times_{ (\rho,a)} \rho||^{(a-1)/2}\times \pi(\phi_-)_1$, où $(\rho,a)$ parcourt l'ensemble des sous-représentations irréductibles
      dans la décomposition de $\phi$ comme ci-dessus, où $b=a$. On constate que la même opération s'applique pour les données
       endoscopiques elliptiques associées à un élément $s$ de $A(\phi)$. 

Cela ramène donc au cas où $\phi$ est trivial sur $\SL(2,\mathbb{C})$. Dans ce cas, $t(\phi)D_{G^*}(\pi(\phi)_1)$ est une 
représentation tempérée irréductible dans le paquet associé à $\phi$. Elle est donc de la forme $\pi(\phi)_\epsilon$ et on
 voit que par les  définitions de ces représentations, ce $\epsilon$ n'est autre que $\epsilon_\phi$. D'où le lemme.\qed

\ 

Pour expliquer un peu ce qui se passe, si $\phi$ est trivial sur $\SL(2,\mathbb{C})$ et n'a pas de multiplicité en tant que
 représentation de $W_F$, alors $\epsilon_\phi$ est trivial car $\pi(\phi)_1$ est une représentation cuspidale. Mais ceci 
 n'est plus vrai en général si $\phi$ est une représentation de $W_F$ avec de la multiplicité car $\pi(\phi)_1$ n'est alors 
 pas une représentation cuspidale et son image par l'involution n'est pas en général elle-même au signe près.

Dans le théorème ci-dessous, les valeurs absolues signifient que l'on prend la représentation en oubliant le signe
 éventuellement introduit par la dualité.

  \begin{thm}
  Soit $\psi$ un paramètre unipotent. Pour tout caractère $\epsilon$ de $A(\psi)$ dont la restriction au centre de 
  $\widehat{G}$ correspond à l'invariant de Hasse, la représentation $\pi(\psi)_\epsilon$ est irréductible et vaut
   $\vert D_G(\pi(\phi)_{\epsilon\epsilon_\phi})\vert$
  où $\phi$ est le morphisme tempéré dual de $\psi$.
  \end{thm}
 \dem  Comme dans \cite{book} paragraphe 7, on identifie le centralisateur de $\phi$ et le centralisateur de  
  $\psi$ (les notations sont celles de l'énoncé). D'après les définitions, on a donc pour tout $\epsilon$ comme dans l'énoncé
  $$
  \epsilon^K(G) D_G(\pi(\phi)_\epsilon)=
  \sum_s \epsilon(s)\,  \mathrm{transfert}\left( D_{H_s} \pi^H_{st}(\phi_s)\right)   ,\eqno(1)
  $$
  où on a noté $H_s$ le groupe de la donnée endoscopique associée à $s$ et $\phi$ et $\phi_s$ la factorisation 
  de $\phi$ par le $L$-groupe de la donnée endoscopique.

   On note, comme dans ce qui précède l'énoncé du théorème, $t(\phi)$ le signe tel que $t(\phi)D_{G^*}(\pi(\phi)_1)$
    est irréductible (ici encore $1$ est le caractère trivial de $A(\phi)$). On sait que $t(\phi)=(-1)^{a_{G^*}-r(\phi)}$ où 
    $r(\phi)$ est le rang du parabolique tel que le module de Jacquet de $\pi(\phi)_1$ pour ce parabolique est formé de représentations cuspidales.

    On note $\sigma(\epsilon)$ le signe tel que 
   $$
   \sigma(\epsilon)t(\phi)\epsilon^K(G) D_G (\pi(\phi)_\epsilon)=:|D(\pi(\phi)_\epsilon)|$$
   soit une représentation irréductible. En revenant à la définition de $\epsilon^K(G)$ et au calcul de signe 
   pour l'involution non normalisée, on vérifie que le signe de Kottwitz remplace la normalisation pour $G$ de 
   l'involution par celle pour $G^*$  et $\sigma(\epsilon)$ est donc $$(-1)^{r(\pi(\phi)_1)-r(\pi(\phi)_\epsilon)}$$ où $r(\pi)$ 
   est le rang du parabolique tel que le module de Jacquet de $\pi$ pour ce parabolique soit cuspidal. 
   Maintenant on peut reprendre mot pour mot \cite{elementaire} 4.2 qui montre que ce signe n'est autre que $\epsilon(s_\phi)$.
   
   On fixe encore $s$ dans le centralisateur de $\phi$ d'où la donnée endoscopique elliptique contenant $s$ 
   et $\phi_s$ comme ci-dessus. On a l'analogue $t(\phi_s)$ de $t(\phi)$. D'après le calcul fait ci-dessus,
    $t(\phi_s) D_{H_s} \pi^H_{st}(\phi_s)$ est la représentation virtuelle  stable associée par Arthur à la factorisation, $\psi_s$ de $\psi$ 
    (où $\psi$ est le morphisme dual de $\phi$) par le groupe endoscopique associé à $s$. 
   Par définition de $\epsilon_\phi$, on a  $t(\phi)=\epsilon_\phi(s)t(\phi_s)$ pour tout $s$. 
   On reporte dans (1) et on trouve
     $$
   |D(\pi(\phi)_\epsilon)|=\sum_s \epsilon(s_\psi)\epsilon_\phi(s)\epsilon(s)\;  \pi^H_{st}(\psi_s).$$
  On utilise encore le fait que $\epsilon_\phi(s_\psi)=1$ pour transformer, pour tout $s\in A(\phi)$
   le signe $\epsilon(s_\psi)\epsilon_\phi(s)\epsilon(s)$ en $(\epsilon\epsilon_\phi)(ss_\psi)$.
    Et on obtient le théorème par les définitions mêmes.
    \qed

   \subsection{Les morphismes unipotents dans le cas archimédiens\label{unipotentarchimedien}}
   Comme \cite{pourhowe} que l'on va utiliser n'est écrit que pour les groupes orthogonaux,
    on se limite à ce cas. Les groupes unitaires peuvent se traiter directement par voie locale comme on le vérifiera dans un article ultérieur.
    
    On fait remarquer aussi que la situation des places complexes est totalement comprise (cf. \cite{GC}) et que seules les places réelles
     nous importent.

   Dans le cas archimédien, un morphisme de $W_F\times \SL(2,\mathbb{C})$ dans le groupe dual de $G$ est dit unipotent 
   s'il est trivial sur le sous-groupe $\mathbb{C}^*$ de $W_F$. Quand on globalise une telle situation comme en
    \cite{pourhowe}, aux places finies on prend des morphismes de $W_{F_v}\times\SL(2,\mathbb{C})\times \SL(2,\mathbb{C})$, 
    triviaux sur la première copie de $\SL(2,\mathbb{C})$ et dont la restriction à $W_{F_v}$ est une somme de caractères quadratiques. 
    Cette situation est  un cas particulier du cas traité dans le paragraphe précédent.
   
   On peut donc ainsi globaliser en connaissant les résultats aux places $p$-adiques  sans l'hypothèse que le groupe 
   est quasi-déployé. Et on connaît, l'existence des paquets d'Arthur avec leur propriété de multiplicité un (cf. \cite{pourhowe}) 
   dans le cas archimédien sous l'hypothèse que le groupe est quasi-déployé. On va donc pouvoir enlever cette hypothèse d'être quasi-déployé.
   
\begin{thm} Soit $\psi$ un morphisme unipotent pour le groupe réel $G$. La représentation $\pi^A_G(\psi)$ restreinte à $G$
 est une somme de représentations irréductibles, somme sans multiplicité.
\end{thm}
   De façon équivalente, le théorème dit que pour tout caractère $\epsilon$ de $A(\psi)$, la représentation $\pi(\psi)_\epsilon$
    est sans multiplicité et disjointe de la représentation $\pi(\psi)_{\epsilon'}$ si $\epsilon'$ est un caractère différent de $A(\psi)$.  
    On fait remarquer la différence avec le cas $p$-adique : $\pi(\psi)_\epsilon$ ici n'est pas irréductible en général, ceci était déjà le 
    cas pour les groupes quasi-déployés. Et cela vient exactement du fait qu'une forme orthogonale n'est pas uniquement
     déterminée par sa dimension, son discriminant et son invariant de Hasse, contrairement au cas $p$-adique.
   
   \dem Pour démontrer le théorème, on globalise la situation en prenant $\mathbb{Q}$ comme corps de nombres.
    Il n'y a donc qu'une seule place archimédienne, celle qui nous intéresse. On globalise $\psi$ de façon
     à avoir en toute place un paramètre unipotent et on connaît donc le théorème en toutes les places sauf la place archimédienne. 
   
   On fixe $\epsilon$ un caractère du centralisateur de $\psi$ tel que $\pi(\psi)_\epsilon$ soit non nul.
    On va démontrer que $\pi(\psi)_\epsilon$ est une somme sans multiplicité de représentations irréductibles.
     Pour cela on utilise la formule globale de \ref{positiviteun} et on sait que $\pi(\psi)_\epsilon$ est une somme de 
     composantes locales de représentations automorphes de carré intégrable. Mais on obtient les représentations 
     automorphes de carré intégrable étudiées en \cite{pourhowe}.  Ces représentations automorphes interviennent avec multiplicité un dans 
     l'espace des représentations automorphes de carré intégrable comme on l'a vérifié en \cite{pourhowe} 3.6. D'où l'assertion.
   
   De plus, on a vérifié en \cite{pourhowe} que deux représentations automorphes de carré intégrable, 
   correspondant à un  paramètre unipotent et qui sont isomorphes en toute place sauf, éventuellement, une place non archimédienne, 
   sont, en fait, isomorphes en toute place. Cela prouve  la fin de l'énoncé comme en {\sl loc.cite}.  Dans l'assertion précédente 
   il est indispensable   que la place où les représentations diffèrent éventuellement soit $p$-adique; aux places archimédiennes
    les représentations     $\pi(\psi)_\epsilon$ ne sont  pas irréductibles même dans le cas unipotent.  \qed

   \end{document}